\pdfoutput=1
\documentclass[a4paper]{article}%
\usepackage{eurosym}
\usepackage{amsfonts}
\usepackage{amssymb}
\usepackage{graphicx}
\usepackage{amsmath}%
\providecommand{\U}[1]{\protect\rule{.1in}{.1in}}
\pdfoutput=1
\newtheorem{theorem}{Theorem}

\newtheorem{conjecture}[theorem]{Conjecture}

\newtheorem{definition}[theorem]{Definition}

\newtheorem{notation}[theorem]{Notation}

\newtheorem{proposition}[theorem]{Proposition}

\newtheorem{summary}[theorem]{Summary}

\ifx\pdfoutput\relax\let\pdfoutput=\undefined\fi
\newcount\msipdfoutput
\ifx\pdfoutput\undefined\else
\ifcase\pdfoutput\else
\ifx\paperwidth\undefined\else
\ifdim\paperheight=0pt\relax\else\pdfpageheight\paperheight\fi
\ifdim\paperwidth=0pt\relax\else\pdfpagewidth\paperwidth\fi
\fi\fi\fi
\begin{document}

\title{On the Contact Numbers of Ball Packings on Various Hexagonal Grids}
\author{Istvan Szalkai\\University of Pannonia, Veszpr\'{e}m, Hungary\thanks{$^{)}$ The author
declares that there is no conflict of interest regarding the publication of
this paper.}$^{)}$\\\texttt{szalkai@almos.uni-pannon.hu}}
\maketitle

\begin{abstract}
We describe the structure of the \textit{different} hexagonal grids in
dimension $d=3$, propose short notation for them, investigate the contact
numbers of ball packings in these grids and share some computational results
up to $200$ balls, using mainly the greedy algorithm. We consider the
octahedral grid, too.

\end{abstract}

\section{\label{Sec-Intro}Introduction}

Considering $n$ balls of the same radius in the $d$ -dimensional Euclidean
space for any but fixed number $n\in\mathbb{N}$ such that each pair of balls
have \textit{at most one} common ("\textit{tangential"}, "\textit{touching}"
or "\textit{kissing}") \textit{point}, the total number of such points is
called the \textbf{contact number} of \textit{this} ball-configuration. For
any $n$ we may ask for configurations which have \textbf{maximal} contact
numbers. These maximal contact numbers are denoted by $c\left(  n,d\right)  $,
we shorten $c\left(  n,3\right)  $ simply by $c\left(  n\right)  $. Many
recent contributions deal with the maximal contact number of balls, using
theoretical and empirical approaches as well (see eg. [AMB11], [B12], [B13],
[BR13], [BK16], [R15], [R16] and [Sz16a]), [BSzSz15] and [H15] deal with
related questions. The possible configurations have many applications e.g. in
material science and other fields of applied physics and chemistry, see e.g.
[BSM11], [HHH12] or [N12].

When packing balls of the same (e.g. \textit{unit}) radius, it is quite
natural to investigate regular packings, eg. on \textit{regular grids}. In the
case $d=2$ (plane) Harboth [H74] proved that \textit{the} (unique)
\textbf{hexagonal grid} (lattice) is the optimal configuration for congruent
circles, see the \textit{blue} circles on Figure 1. This means that the
centers of the circles fit on the grid (a lattice in fact)
\begin{equation}
\left\{  i\cdot\left[  2,0\right]  ^{T}+j\cdot\left[  1,\sqrt{3}\right]
^{T}~:\text{\quad}i,j\in\mathbb{Z}\right\}  ~\text{.} \label{2dim}%
\end{equation}

In dimension $d=3$ the hexagonal grids can be obtained as follows. First,
place some balls in \textit{planes} in the above planar method, we call these
sets \textbf{layers.} Second, put such layers onto each other: the layers are
translated with a vector like $\ell_{1}:=\left[  1,\sqrt{1/3},\sqrt
{8/3}\right]  ^{T}$ to get the neighbour layer (like in the usual placement of
melons), as the \textit{blue} and \textit{red} circles show in Figure 1 (top view).

\textit{However}, when considering three or more layers, we have several
\textit{different} possibilities. Think the blue layer in the \textit{middle}
and the red one \textit{above} it. Then, for the layer \textit{below} the blue
circles we have two \textit{different} possibilities. \textit{Either} just put
the balls "exactly below" the red ones, i.e. translating the red layer with
the vector $t:=\left[  0,0,-2\cdot\sqrt{8/3}\right]  ^{T}$, in which case both
the layers above and below the blue one are red. The \textit{other}
possibility is to move the \textit{blue} layer with the vector $-\ell_{1}$
resulting the green circles. Equivalently, rotating the red layer by $60^{o}$
and translate it with the vector $t$ we get the green layer.

\begin{center}%
\raisebox{-0cm}{\parbox[b]{11.1852cm}{\begin{center}
\ifcase\msipdfoutput
\includegraphics[
natheight=16.757999cm,
natwidth=22.225401cm,
height=8.4504cm,
width=11.1852cm
]%
{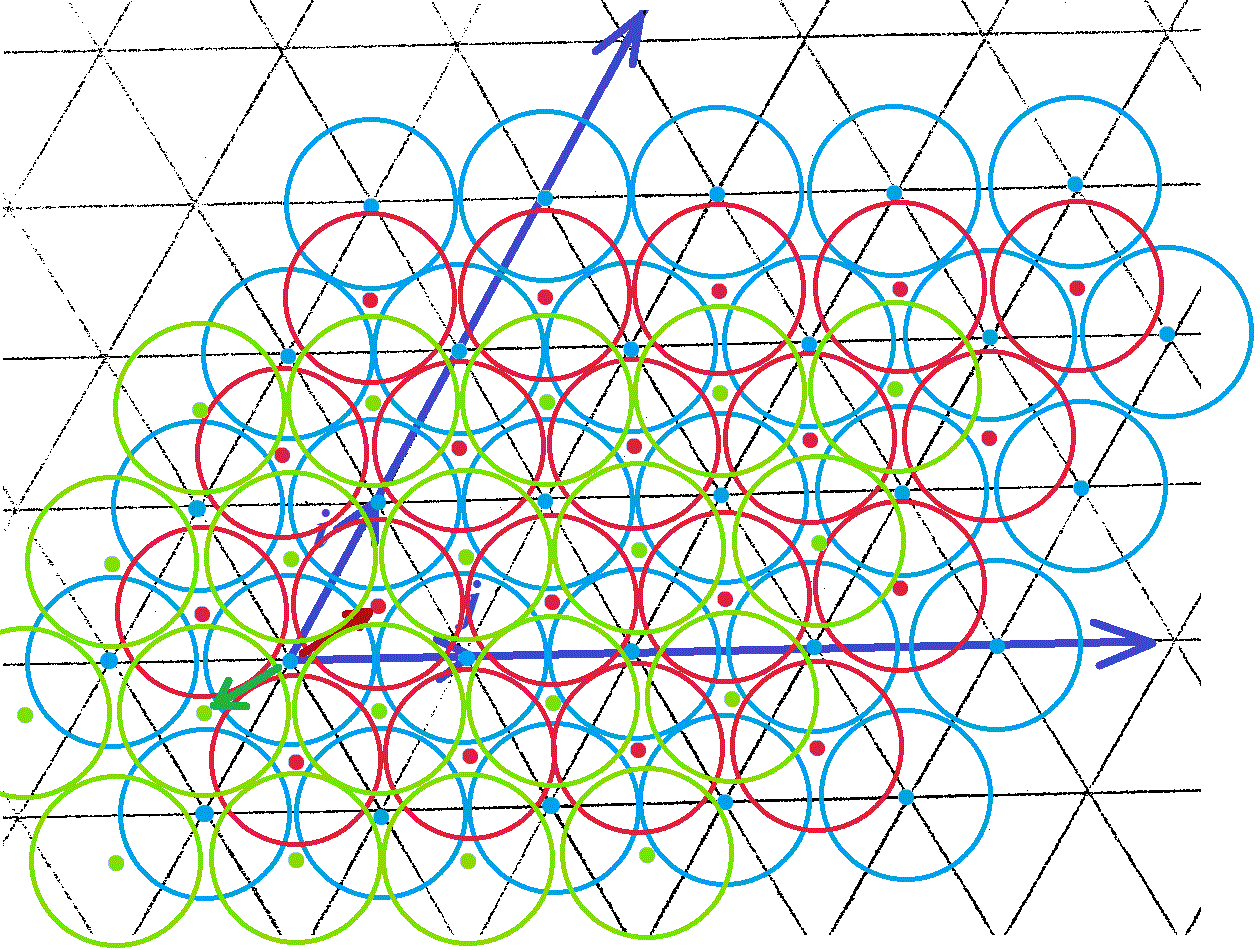}%
\else
\includegraphics[
natheight=8.450400cm,
natwidth=11.185200cm,
height=8.4504cm,
width=11.1852cm
]%
{3szog-racs-korok-kekpiroszold.png}%
\fi
\\
\textbf{Figure 1:} \textit{Possible neighbouring layers (top view)}%
\end{center}}}

\end{center}

We propose to choose and fix a \textit{middle} (blue) layer for forthcoming
constructions, since when packing the next ball (step by step) we can place it
both above and below the configuration we already have.

In general, when placing several layers \textit{above each other}, for each
next layer we can chose the vector either $\ell_{1}=L_{v}+L_{h}$ or $\ell
_{2}=L_{v}-L_{h}$ to move the previous layer to get the next one \textit{above
it}, where
\begin{equation}
L_{v}:=\left[  0,0,\sqrt{8/3}\right]  ^{T}\text{~,\quad}L_{h}:=\left[
1,\sqrt{1/3},0\right]  ^{T} \label{Lv_es_Lh}%
\end{equation}
are the \textbf{vertical} and \textbf{horizontal} translating vectors.
Similarly, if we plan a next layer \textit{below the bottom layer} we may
choose either $-\ell_{1}=-L_{v}-L_{h}$ or $-\ell_{2}=-L_{v}+L_{h}$ to move the
bottom layer to get the next one \textit{below it}. A general notation for all
the obtainable grids is discussed in Section \ref{Sec-Notations}. The usual
\textit{hexagonal lattice} can be obtained when taking always the same vector,
e.g. $\ell_{1}$ between all consecutive layers.

Though the contacts of balls depend only on neighbouring layers, but when we
have to place \textit{fixed number} ($n$) of balls, different grids
\textit{often} provide different configurations and different maximal contact
numbers. Examples in [Sz16b] show even $30\%$ differences in different grids!
This problem and examples are discussed in Section \ref{Sec-CompRes} and in [Sz16b].

An important question is whether $c\left(  n\right)  $ can be achieved for all
$n\in\mathbb{N}$ in one of the hexagonal grids. Since all the grids are
$12$-regular (each ball has $12$ touching neighbours) we suspect that the
answer is \textit{yes}: $c_{grids}\left(  n\right)  =6n-o\left(  n\right)  $
might hold for $n\rightarrow\infty$\textit{.} For example, none of the
contructions shown in [AMB11, p.35] and in [BK16] for $c\left(  6\right)  =12$
do exist in any grid, but the configuration shown in Figure 2 (top and
perspective view) yield the same contact number.

\begin{center}%
\raisebox{-0cm}{\ifcase\msipdfoutput
\includegraphics[
natheight=7.938600cm,
natwidth=6.879800cm,
height=4.0396cm,
width=3.5102cm
]%
{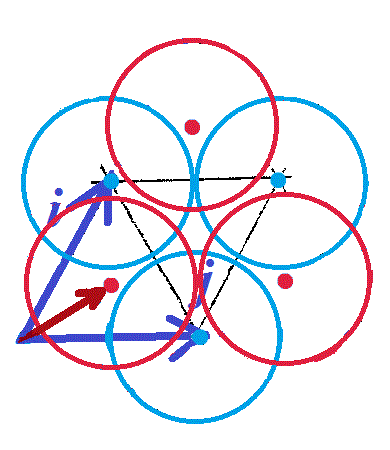}%
\else
\includegraphics[
natheight=4.039600cm,
natwidth=3.510200cm,
height=4.0396cm,
width=3.5102cm
]%
{Gn0601-kekpiros.png}%
\fi
}
\qquad%
\raisebox{-0cm}{\ifcase\msipdfoutput
\includegraphics[
natheight=5.115900cm,
natwidth=6.701900cm,
height=3.9078cm,
width=5.1006cm
]%
{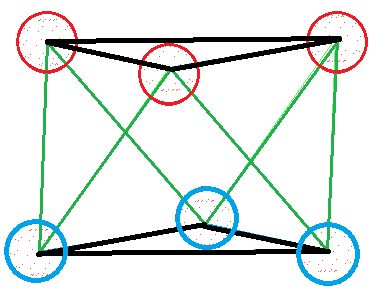}%
\else
\includegraphics[
natheight=3.907800cm,
natwidth=5.100600cm,
height=3.9078cm,
width=5.1006cm
]%
{Gn0601-jav2-.png}%
\fi
}
\medskip

\textbf{Figure 2:} $c\left(  6\right)  =12$ \textit{on a grid }\bigskip
\end{center}

\noindent Exact values of $c\left(  n\right)  $ for $n\leq20$ can be found
e.g. in [AMB11, p.38] and [BK16, p.5] and some lower bounds for $n\leq27$ in
[Sz16a]. Our present computations (see Section \ref{Sec-CompRes} and [Sz16b])
show that for all these $n$ the (maximal) value of $c\left(  n\right)  $ can
be achieved in some \textit{appropriate} grid, possibly not by the greedy
algorithm. Examples for this phenomane are the cases $n=14$ and $n=15$~:\ the
greedy algorithm (see Table 2) found less contacts, but other algorithm (see
Section \ref{Sec-CompRes}) found the exact values:

\begin{center}%
\raisebox{-0cm}{\ifcase\msipdfoutput
\includegraphics[
natheight=8.819400cm,
natwidth=11.288500cm,
height=5.1533cm,
width=5.7859cm
]%
{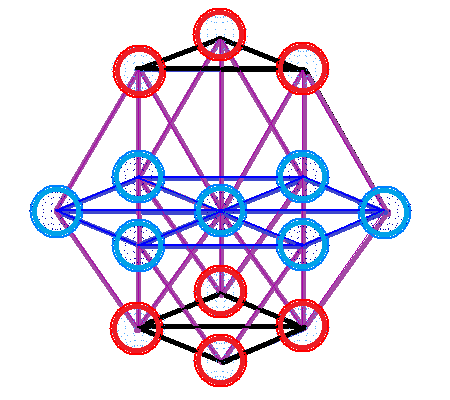}%
\else
\includegraphics[
natheight=5.153300cm,
natwidth=5.785900cm,
height=5.1533cm,
width=5.7859cm
]%
{G14H3g-H1d.png}%
\fi
}
\raisebox{-0cm}{\ifcase\msipdfoutput
\includegraphics[
natheight=7.055600cm,
natwidth=11.288500cm,
height=5.1533cm,
width=5.7859cm
]%
{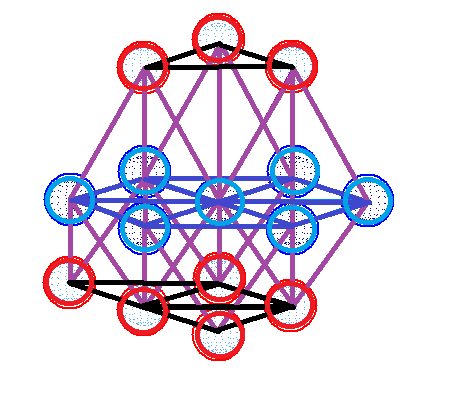}%
\else
\includegraphics[
natheight=5.153300cm,
natwidth=5.785900cm,
height=5.1533cm,
width=5.7859cm
]%
{G15H3G4-H13d.png}%
\fi
}

\textbf{Figure 3:} $c\left(  14\right)  =40$ and $c\left(  15\right)  =44$
\textit{on a grid }\bigskip
\end{center}

So we might state the following:

\begin{conjecture}
\label{Con-Hex-Sk}The exact (maximal) value for each $n\in\mathbb{N}$ can be
achieved in some hexagonal grid. $\Box$
\end{conjecture}

We have to emphasize that different $n$ numbers often require different grids
for achieving the maximal contact number. Our computational results for
$n\leq200$ are explained in Section \ref{Sec-CompRes}, the details can be
found in [Sz16b].

We consider also the octahedral lattice at the end of the forthcoming Sections
and in [Sz16b], comparing to [B12] and [R16].

\section{\label{Sec-Notations}Notations}

\begin{notation}
For \emph{any} geometrical point $P\in\mathbb{R}^{3}$ we denote the usual
Cartesian coordinates of $P$ by $\left[  x,y,z\right]  _{D}^{T}$ or $\left[
P\right]  _{D}$~, we often leave the subscript $D$ (and $T$ is for transposing
the vector). \newline If $P\in\mathcal{G}_{\overrightarrow{\varepsilon}}$ is
an element of a \emph{hexagonal} grid $\mathcal{G}_{\overrightarrow
{\varepsilon}}$~, determined be the vector $\overrightarrow{\varepsilon}$ (see
Definition \ref{Def-Geps}), the \emph{inner} hexagonal coordinates are denoted
by $\left[  P\right]  _{\mathcal{G}}$ or $\left[  P\right]  _{\overrightarrow
{\varepsilon}}$~, or simply by $\left[  P\right]  _{H}$ when $\overrightarrow
{\varepsilon}$ is clear from the context.\newline In each grid we
have\quad$\left[  0,0,0\right]  _{\mathcal{G}}^{T}=\left[  0,0,0\right]
_{D}^{T}$~.\quad$\Box$
\end{notation}

\begin{definition}
\label{Def-Geps}For any infinite (or finite) sequence \newline$\overrightarrow
{\varepsilon}=\left(  ...,\varepsilon_{-m},...,\varepsilon_{-1},\varepsilon
,\varepsilon_{1},...,\varepsilon_{n},...\right)  $ such that $\varepsilon
_{0}=0$ \ and $\varepsilon_{k}\in\left\{  -1,+1\right\}  $ for $k\in
\mathbb{Z}\backslash\left\{  0\right\}  $ we define the \textbf{hexagonal
grid} $\mathcal{G}_{\overrightarrow{\varepsilon}}\subset\mathbb{R}^{3}$ as
follows.\newline The $0$ -st layer (as in (\ref{2dim})) is
\begin{equation}
\mathcal{L}_{0}:=\left\{  i\cdot\left[
\begin{array}
[c]{c}%
2\\
0\\
0
\end{array}
\right]  _{D}+j\cdot\left[
\begin{array}
[c]{c}%
1\\
\sqrt{3}\\
0
\end{array}
\right]  _{D}~:~i,j\in\mathbb{Z}\right\}
\end{equation}
and for $P\in\mathcal{L}_{0}$ the hexagonal coordinates are $\left[  P\right]
_{\mathcal{G}}^{T}:=\left[  i,j,0\right]  ^{T}$~.\newline For $k>0$
($k\in\mathbb{Z}$) we define the $k$ -th layer
\begin{equation}
\mathcal{L}_{k}:=\left\{  R+L_{v}+\varepsilon_{k}L_{h}~:~R\in\mathcal{L}%
_{k-1}\right\}
\end{equation}
(see (\ref{Lv_es_Lh}) for $L_{v}$ and $L_{h}$), the hexagonal coordinates for
$P=R+L_{v}+\varepsilon_{k}L_{h}$ and $\left[  R\right]  _{\mathcal{G}}%
^{T}:=\left[  i,j,k-1\right]  ^{T}$ are $\left[  P\right]  _{\mathcal{G}}%
^{T}:=\left[  i,j,k\right]  ^{T}$~.\newline For $k<0$ ($k\in\mathbb{Z}$) we
let the $k$ -th layer
\begin{equation}
\mathcal{L}_{k}:=\left\{  R-L_{v}+\varepsilon_{k}L_{h}~:~R\in\mathcal{L}%
_{k+1}\right\}
\end{equation}
and the hexagonal coordinates for $P=R-L_{v}+\varepsilon_{k}L_{h}$ and
$\left[  R\right]  _{\mathcal{G}}^{T}:=\left[  i,j,k+1\right]  ^{T}$ are
$\left[  P\right]  _{\mathcal{G}}^{T}:=\left[  i,j,k\right]  ^{T}$~.\newline
Finally we put
\begin{equation}
\mathcal{G}_{\overrightarrow{\varepsilon}}:=\bigcup_{k\in\mathbb{Z}%
}\mathcal{L}_{k}\text{\ .\quad}\Box\label{G_eps=}%
\end{equation}

\end{definition}

Clearly in $\left[  i,j,k\right]  _{\mathcal{G}}^{T}$\quad$k$ denotes the
layer number, i.e. "\textit{vertical}" coordinate, while $i,j$ are
"\textit{horizontal}" ones as in (\ref{2dim}). The relation whether
"\textit{any two balls with centers in hexagonal coordinates} $\left[
i_{1},j_{1},k_{1}\right]  _{\mathcal{G}}^{T}$ \textit{and} $\left[
i_{2},j_{2},k_{2}\right]  _{\mathcal{G}}^{T}$ \textit{contact each other or
not}" can be easily decided.

The usual \textbf{hexagonal lattice} can be obtained when taking always the
same direction, e.g. for $\overrightarrow{\varepsilon}_{HL}=\left(
...,-1,0,+1,...\right)  $. This lattice is
\begin{equation}
\mathcal{G}_{HEX}=\left\{  x\cdot\overrightarrow{i}+y\cdot\overrightarrow
{i}+z\cdot\overrightarrow{\ell_{1}}:x,y,z\in\mathbb{Z}\right\}  \label{G_HEX=}%
\end{equation}
where
\begin{equation}
\overrightarrow{i}=\left[
\begin{array}
[c]{c}%
2\\
0\\
0
\end{array}
\right]  _{D}\text{ ,\quad}\overrightarrow{j}=\left[
\begin{array}
[c]{c}%
1\\
\sqrt{3}\\
0
\end{array}
\right]  _{D}\text{ ,\quad}\ell_{1}=\left[
\begin{array}
[c]{c}%
1\\
\sqrt{1/3}\\
\sqrt{8/3}%
\end{array}
\right]  \text{ .}%
\end{equation}

In the case our balls all are contained in $\bigcup\limits_{k=t_{1}}^{t_{2}%
}\mathcal{L}_{k}$ for some fixed $t_{1},t_{2}\in\mathbb{Z}$, then only the
values $\varepsilon_{t_{1}},...,\varepsilon_{t_{2}}$ are interesting for us,
so $\overrightarrow{\varepsilon}$ can be chosen a finite sequence and can be
shortened in a single (rational) number.

\begin{notation}
For finite sequences $\overrightarrow{\varepsilon}=\left(  \varepsilon_{t_{1}%
},...,\varepsilon_{t_{2}}\right)  $ \ for which $\varepsilon_{k}\in\left\{
-1,+1\right\}  $, for any $t_{1},t_{2}\in\mathbb{Z}$ we write
\begin{equation}
type\left(  \overrightarrow{\varepsilon}\right)  :=\sum_{k=t_{1}}^{t_{2}%
}2^{s\left(  \varepsilon_{k}\right)  }\text{\quad where }s\left(
\varepsilon_{k}\right)  =\left\{
\begin{array}
[c]{l}%
0\text{\quad if }\varepsilon_{k}=-1\\
1\text{\quad if }\varepsilon_{k}=+1
\end{array}
\right.  \text{ .\quad}\Box
\end{equation}

\end{notation}

\begin{proposition}
For any $k\in\mathbb{Z}$ and $P\in\mathcal{G}$ , $P=\left[  i,j,k\right]
_{\overrightarrow{\varepsilon}}^{T}\in\mathcal{L}_{k}$ we have \medskip
\newline\textbf{(i)} the \emph{Hexagonal} coordinates of $P$ are: for
$k\geq0$
\begin{equation}
P\left[  i,j,k\right]  _{\overrightarrow{\varepsilon}}^{T}=P\left[
i,j,0\right]  ^{T}+kL_{v}+L_{h}\sum_{t=1}^{k}\varepsilon_{t}%
\end{equation}
(where $\sum\limits_{t=1}^{0}\varepsilon_{t}=0$), and for $k<0$
\begin{equation}
P\left[  i,j,k\right]  _{\overrightarrow{\varepsilon}}^{T}=P\left[
i,j,0\right]  ^{T}+kL_{v}+L_{h}\sum_{t=k}^{-1}\varepsilon_{t}\text{ .}%
\end{equation}
\medskip\textbf{(ii) }the \emph{Cartesian} coordinates of the point $P$ are
for each $k\in\mathbb{Z}$%
\begin{equation}
\left[  P\right]  _{D}=i\cdot\left[
\begin{array}
[c]{c}%
2\\
0\\
0
\end{array}
\right]  _{D}+j\cdot\left[
\begin{array}
[c]{c}%
1\\
\sqrt{3}\\
0
\end{array}
\right]  _{D}+k\cdot\left[
\begin{array}
[c]{c}%
0\\
0\\
\sqrt{8/3}%
\end{array}
\right]  _{D}+S_{k}\cdot\left[
\begin{array}
[c]{c}%
1\\
\sqrt{1/3}\\
0
\end{array}
\right]  _{D} \label{PD=}%
\end{equation}
where
\begin{equation}
S_{k}=\left\{
\begin{array}
[c]{l}%
\sum_{t=1}^{k}\varepsilon_{t}\\
0\\
\sum_{t=k}^{-1}\varepsilon_{t}%
\end{array}
\right.  ~\left.
\begin{array}
[c]{c}%
\text{for }k>0\\
\text{for }k=0\\
\text{for }k<0
\end{array}
\right.  \text{~.\quad}\Box
\end{equation}

\end{proposition}

\noindent Similar (also triangular) \textit{barycentric} coordinate system is
also used in chemistry, see e.g. [Sz99].

The \textit{octahedral} grid is, in fact, a lattice and is unique. Any layer
is a lattice of squares of side $2$, the vector $\left[  1,1,\sqrt{2}\right]
^{T}$ moves layers to the next one, i.e. we have:

\begin{definition}
The \textbf{octahedral lattice} is
\begin{equation}
\mathcal{G}_{OCT}:=\left\{  x\cdot\left[
\begin{array}
[c]{c}%
2\\
0\\
0
\end{array}
\right]  +y\cdot\left[
\begin{array}
[c]{c}%
0\\
2\\
0
\end{array}
\right]  +z\cdot\left[
\begin{array}
[c]{c}%
1\\
1\\
\sqrt{2}%
\end{array}
\right]  ~:~x,y,z\in\mathbb{Z}\right\}  \text{ . }\Box\label{G_OCT=}%
\end{equation}

\end{definition}

\section{\label{Sec-Former}Former results}

Exact values of $c\left(  n\right)  $ for $n\leq19$ can be found e.g. in
[AMB11, p.38] and [BK16, p.5], some \textit{lower bounds} for $20\leq n\leq27$
in [Sz16a]:

\begin{center}%
\begin{tabular}
[c]{|r||r|r|r|r|r|r|r|r|r|r|}\hline
$n$ & $0$ & $1$ & $2$ & $3$ & $4$ & $5$ & $6$ & $7$ & $8$ & $9$\\\hline\hline
$0$ & $-$ & $0$ & $1$ & $3$ & $6$ & $9$ & $12$ & $15$ & $18$ & $21$\\\hline
$10$ & $25$ & $29$ & $33$ & $36$ & $40$ & $44$ & $48$ & $52$ & $56$ &
$60$\\\hline
$20$ & $64$ & $67$ & $72$ & $76$ & $80$ & $84$ & $87$ & $90$ & $-$ &
$-$\\\hline
\end{tabular}
\medskip

\textbf{Table 1} \textit{Results from }[BK16] ($1\leq n\leq19$) \textit{and
from} [Sz16a] ($20\leq n\leq27$) \bigskip
\end{center}

\noindent Using a construction in the \textit{octahedral} lattice [B12],
[BR13] and [BK16] states
\begin{equation}
c\left(  n\right)  \geq4k^{3}-6k^{2}+2k=d\left(  k\right)  \thickapprox
6n-o\left(  n^{2/3}\right)  \label{BK16T41}%
\end{equation}
for $n=\dfrac{2k^{3}+k}{3}$ ($k\in\mathbb{N}$), that is

\begin{center}%
\begin{tabular}
[c]{|r||r|r|r|r|r|r|r|}\hline
$k$ & $1$ & $2$ & $3$ & $4$ & $5$ & $6$ & $...$\\\hline\hline
$\mathbf{n}$ & $\mathbf{1}$ & $\mathbf{6}$ & $\mathbf{19}$ & $\mathbf{44}$ &
$\mathbf{85}$ & $\mathbf{146}$ & $\mathbf{...}$\\\hline
$c\left(  n\right)  \geq$ & $0$ & $12$ & $60$ & $168$ & $360$ & $660$ &
$...$\\\hline
\end{tabular}
\medskip

\textbf{Table 2}\textit{\ Summary of }[B12]\bigskip
\end{center}

\noindent\lbrack R16] contains a figure for $k=4$ and some values for these
$n$ and their $c\left(  n\right)  $ are shown in Table 2.

Summary \ref{Osszegezes} gives a comparison of the results obtained by
different methods.

\section{\label{Sec-CompRes}Computational results}

In [Sz16a] we forced the computer to check \textit{all the cases} in a
$3\times3\times3$ grid, i.e. $3$ layers and $-1\leq i,j\leq+1$ for all
$n\leq27$, running time varied from some seconds to $1-2$ hours. Since this
kind of total checkings in $4\times4\times4$ required several days running
time, we turned to the \textit{greedy algorithm}, which terminated in some
seconds even for $n>200$.

In our recent computations we investigated $n$ many balls up to $n\leq200$ in
the grids $\mathcal{G}_{i}$ listed in [Sz16b] in file
\textsc{GnMohoH3e-160626-0050.*} and defined in Section \ref{Sec-Notations}.
Our algorithm had the possibility to use $31$ layers, starting at layer $0$,
but all the balls which have been finally chosen plus their neighbours
(possible further ones for larger $n$ configurations) occupied only the layers
$\mathcal{L}_{-4},...,\mathcal{L}_{4}$~. In other words, we had to investigate
the sequences $\overrightarrow{\varepsilon_{0}}=\left(  \varepsilon
_{-4},...,\varepsilon_{-1},0,\varepsilon_{1},...\varepsilon_{4}\right)  $
only, which means $2^{8-1}=128$ cases, since, by symmetry we assumed
$\varepsilon_{1}=+1$~.

Greedy algorithm means, that we started with an arbitrary ball, and in each
step we chosed the next one which has the most neighbour (touching ball) among
the old ones (we have chosen in previous steps). The greedy method implies
that from the configuration for $200$ balls we can reconstruct each
configurations for each $n\leq200$ - we have to consider the balls
$B_{1},...,B_{n}$ only. So, one can easily read out each configuration in each
grid for each $n\leq200$ from the files in [Sz16b]. Table 3 contains the
maximum contact number for each $n$, we have investigated all the $128$ grids.
\medskip

\begin{center}%
\begin{tabular}
[c]{|r||r|r|r|r|r|r|r|r|r|r|}\hline
$c\left(  n\right)  $ & $\mathbf{0}$ & $\mathbf{1}$ & $\mathbf{2}$ &
$\mathbf{3}$ & $\mathbf{4}$ & $\mathbf{5}$ & $\mathbf{6}$ & $\mathbf{7}$ &
$\mathbf{8}$ & $\mathbf{9}$\\\hline\hline
$\mathbf{0}$ & $-$ & $0$ & $1$ & $3$ & $6$ & $9$ & $\mathbf{11}$ & $15$ & $18$
& $21$\\\hline
$\mathbf{10}$ & $25$ & $29$ & $33$ & $36$ & $\mathbf{39}$ & $\mathbf{43}$ &
$48$ & $52$ & $56$ & $60$\\\hline
$\mathbf{20}$ & $64$ & $68$ & $72$ & $75$ & $79$ & $84$ & $89$ & $93$ & $97$ &
$102$\\\hline
$\mathbf{30}$ & $106$ & $110$ & $114$ & $119$ & $123$ & $126$ & $130$ & $135$
& $140$ & $145$\\\hline
$\mathbf{40}$ & $150$ & $153$ & $157$ & $162$ & $167$ & $172$ & $177$ & $183$
& $187$ & $191$\\\hline
$\mathbf{50}$ & $195$ & $200$ & $205$ & $210$ & $214$ & $218$ & $222$ & $227$
& $232$ & $236$\\\hline
$\mathbf{60}$ & $242$ & $247$ & $251$ & $257$ & $261$ & $265$ & $271$ & $275$
& $280$ & $284$\\\hline
$\mathbf{70}$ & $288$ & $293$ & $298$ & $303$ & $308$ & $312$ & $317$ & $322$
& $328$ & $332$\\\hline
$\mathbf{80}$ & $337$ & $342$ & $348$ & $352$ & $356$ & $360$ & $365$ & $369$
& $375$ & $380$\\\hline
$\mathbf{90}$ & $385$ & $389$ & $394$ & $398$ & $403$ & $408$ & $414$ & $419$
& $424$ & $428$\\\hline
$\mathbf{100}$ & $433$ & $438$ & $444$ & $448$ & $453$ & $458$ & $463$ & $468$
& $473$ & $477$\\\hline
$\mathbf{110}$ & $481$ & $487$ & $491$ & $496$ & $501$ & $505$ & $510$ & $514$
& $519$ & $524$\\\hline
$\mathbf{120}$ & $530$ & $535$ & $541$ & $546$ & $551$ & $555$ & $559$ & $563$
& $568$ & $573$\\\hline
$\mathbf{130}$ & $578$ & $583$ & $589$ & $594$ & $600$ & $605$ & $610$ & $615$
& $620$ & $625$\\\hline
$\mathbf{140}$ & $630$ & $633$ & $638$ & $643$ & $648$ & $652$ & $658$ & $663$
& $669$ & $674$\\\hline
$\mathbf{150}$ & $679$ & $684$ & $690$ & $695$ & $701$ & $706$ & $712$ & $717$
& $723$ & $727$\\\hline
$\mathbf{160}$ & $731$ & $735$ & $741$ & $745$ & $750$ & $755$ & $760$ & $766$
& $771$ & $777$\\\hline
$\mathbf{170}$ & $782$ & $788$ & $793$ & $798$ & $802$ & $806$ & $810$ & $815$
& $820$ & $825$\\\hline
$\mathbf{180}$ & $830$ & $834$ & $839$ & $844$ & $849$ & $855$ & $860$ & $865$
& $871$ & $877$\\\hline
$\mathbf{190}$ & $882$ & $888$ & $894$ & $899$ & $905$ & $909$ & $914$ & $920$
& $925$ & $931$\\\hline
$\mathbf{200}$ & $935$ &  &  &  &  &  &  &  &  & \\\hline
\end{tabular}
\medskip

\textbf{Table 3} \textit{Hexagonal lower bounds by greedy algorithm}\bigskip
\end{center}

We also made a run in the \textit{octahedral} lattice, details can be found in
[Sz16b] in files \textsc{MH5b-200b-en.*}. txt Table 4 shows the cases when
octahedral greedy algorithm resulted better bound than the hexagonal one in
Table 3: \medskip

\begin{center}
$%
\begin{tabular}
[c]{|r|r|r||r|r|r||r|r|r|r|r|r|r|}\hline
$n$ & $\mathbf{14}$ & $\mathbf{15}$ & $\mathbf{57}$ & $\mathbf{58}$ &
$\mathbf{59}$ & $\mathbf{176}$ & $\mathbf{177}$ & $\mathbf{178}$ &
$\mathbf{179}$ & $\mathbf{180}$ & $\mathbf{181}$ & $\mathbf{182}$\\\hline
\textbf{Hexa} & $39$ & $43$ & $227$ & $232$ & $236$ & $810$ & $815$ & $820$ &
$825$ & $830$ & $834$ & $839$\\\hline
\textbf{Octa} & $\mathbf{40}$ & $\mathbf{44}$ & $\mathbf{228}$ &
$\mathbf{233}$ & $\mathbf{237}$ & $\mathbf{811}$ & $\mathbf{817}$ &
$\mathbf{822}$ & $\mathbf{828}$ & $\mathbf{833}$ & $\mathbf{837}$ &
$\mathbf{841}$\\\hline
\end{tabular}
\medskip$

\textbf{Table 4}\textit{\ The cases when octahedral is better than hexagonal}
\bigskip
\end{center}

We can shortly summarize our experiments as:

\begin{summary}
\label{Osszegezes}Based on \textbf{Table 3} (i.e. the greedy algorithm in
hexagonal grids), we can state the followings. \newline(i) for $n\leq19$ Table
3 gives the so far known results except for $n=6$ , but this case is cured in
Figure 2,\newline(ii) for $n\leq200$ Table 3 gives the best lower bounds with
the following exceptions: \newline- for $n=6$ see Figure 2 (so $12\leq
c\left(  6\right)  $),\newline- for $n=14$, $15$ see Figure 3 and Table 4 (so
$40\leq c\left(  14\right)  $ and $44\leq c\left(  15\right)  $),\newline- for
$n=23$, $24$ see Table 1 and [Sz16a] (so $76\leq c\left(  23\right)  $ and
$80\leq c\left(  24\right)  $),\newline- for $n=44$ see Table 2 (so $168\leq
c\left(  44\right)  $),\newline- for $n=57-59$ see Table 4 (so $228\leq
c\left(  57\right)  $, $233\leq c\left(  58\right)  $ and $237\leq c\left(
59\right)  $),\newline- for $n=146$ see Table 2 (so $660\leq c\left(
146\right)  $),\newline- for $n=176-182$ see Table 4 \newline(so $811\leq
c\left(  176\right)  $, $817\leq c\left(  177\right)  $, $822\leq c\left(
178\right)  $, $828\leq c\left(  179\right)  $, $833\leq c\left(  180\right)
$, $837\leq c\left(  181\right)  $ and $841\leq c\left(  182\right)  $).
$\Box$
\end{summary}

Though for the exceptional cases $n\geq44$ we do not have constructions in
hexagonal grids at this moment, as good as in Tables 2 and 4, we still believe
in Conjecture \ref{Con-Hex-Sk}.

Further experiments are in progress.

\bigskip

\section*{References}

.

\textbf{[AMB11]} \textbf{Arkus,N., Manoharan,V., Brenner,M.:} \textit{Deriving
Finite Sphere Packings}, SIAM J. Discrete Math. Vol. 25, No. 4, pp. 1860-1901,
\newline\textrm{http://arxiv.org/abs/1011.5412v2}~.\medskip

\textbf{[B12]} \textbf{Bezdek,K.:} \textit{Contact Numbers for Congruent
Sphere Packings in Euclidean 3-Space}, Discrete Comput Geom (2012)
48:298--309, DOI 10.1007/s00454-012-9405-9\medskip

\textbf{[B13]} \textbf{Bezdek,K.:} \textit{Lectures on Sphere Arrangements -
the Discrete Geometric Side}, Springer, 2013.\medskip

\textbf{[BR13]} \textbf{Bezdek,K., Reid,S.:} \textit{Contact graphs of unit
sphere packings revisited}, J. Geom. 104 (2013), no. 1, 57-83.\medskip

\textbf{[BSzSz15]} \textbf{Bezdek,K., Szalkai,B., Szalkai,I.:} \textit{On
contact numbers of totally separable unit sphere packings}, Discrete
Mathematics, Vol. 339, No. 2 (November, 2015), pp. 668-676.\medskip

\textbf{[BK16]} \textbf{Bezdek,K., Khan,M.A.:} \textit{Contact numbers for
sphere packings}, \textrm{http://arxiv.org/abs/1601.00145}~.\medskip

\textbf{[BSM11]} \textbf{Blair,D., Santangelo,C.D., Machta,J.:}
\textit{Packing Squares in a Torus}, J. Statistical Mechanics: Theory and
Experiment, P01018 (2012),
\textrm{http://iopscience.iop.org/1742-5468/2012/01/P01018/article},
\newline\textrm{http://arxiv.org/abs/1110.5348}~.\medskip

\textbf{[H74]} \textbf{Harborth, H.:} \textit{L\"{o}sung zu Problem} 664A,
Elem. Math. 29, 14--15 (1974).\medskip

\textbf{[H15]} \textbf{Holmes-Cerfon,M.:} \textit{Enumerating nonlinearly
rigid sphere packings}, \textrm{http://arxiv.org/abs/1407.3285v2}~.\medskip

\textbf{[HHH12]} \textbf{Hoy,R.S., Harwayne-Gidansky,J., O'Hern,C.S.:}
\textit{Structure of finite sphere packings via exact enumeration:
implications for colloidal crystal nucleation}, Physical Review E 85 (2012),
051403.\medskip

\textbf{[N12]} \textbf{Neves, E.J.:} \textit{A discrete variational problem
related to Ising droplets at low temperatures}, J. Statistical Physics, 80
(1995), 103-123.\medskip

\textbf{[R15]} \textbf{Reid, S.: }\textit{Regular Totally Separable Sphere
Packings}, RESEARCH\textperiodcentered June 2015, DOI:
10.13140/RG.2.1.2771.2161, \textrm{http://arxiv.org/abs/1506.04171}~.\medskip

\textbf{[R16]} \textbf{Reid, S.:} \textit{On Contact Numbers of Finite Lattice
Sphere Packings and the Maximal Coordination of Monatomic Crystals},
\textrm{http://arxiv.org/pdf/1602.04246.pdf}~.\medskip

\textbf{[Sz99]} \textbf{Szalkai,I.:} \textit{Handling Multicomponent Systems
in} $\mathbb{R}^{n}$, I., J. Math.Chem. 25\ (1999), 31-46.\emph{\medskip}

\textbf{[Sz16a]} \textbf{Szalkai,I.:} \textit{On Contact Numbers of Finite
Lattice Sphere Packings of 20-27 Balls},
\textrm{http://arxiv.org/abs/1603.05390}~.\emph{\medskip}

\textbf{[Sz16b] Szalkai,I.:} \textit{Computational results in }$3$%
\textit{\ dimensional hexagonal grids - datasets,}
\textrm{http://math.uni-pannon.hu/\symbol{126}szalkai/GnMH-160630.zip}%
~.\medskip
\end{document}